\documentclass[aop,MSNbibl,nameyear,seceqn,dvips]{arximspdf}

\doi{10.1214/12-AOP787} 
\volume{41}
\issue{5}
\pubyear{2013}
\firstpage{3494}
\lastpage{3517}

\makeatletter

\newcommand{\eqref}[1]{(\ref{#1})}
\newtheorem{lemma}{Lemma}[section]
\newtheorem{Lemma}{Lemma}[section]

\newtheorem{theorem}[Lemma]{Theorem}
\newproclaim{defn}[prop]{Definition}

\newtheorem{proposition}[Lemma]{Proposition}
\newproclaim{remark}[Lemma]{Remark}
\newproclaim{example}[Lemma]{Example}
\newproclaim{ass}{Assumption}

\newcommand{\reals}{\mathbb{R}}
\newcommand{\bbr}{\reals}
\newcommand{\vep}{\varepsilon}

\newcommand{\BX}{\mathbf{X}}

\newcommand{\one}{\mathbf{1}}

\newcommand{\eid}{\stackrel{d}{=}}

\makeatother

\begin{document}
\begin{frontmatter}

\title{Is the location of the supremum of a stationary
process nearly uniformly distributed?\thanksref{T1}}
\runtitle{Location of the supremum}

\thankstext{T1}{Supported in part by the ARO
Grant W911NF-07-1-0078, NSF Grant DMS-10-05903 and NSA Grant
H98230-11-1-0154 at Cornell University.}

\begin{aug}
\author[A]{\fnms{Gennady} \snm{Samorodnitsky}\corref{}\ead[label=e1]{gs18@cornell.edu}}
\and
\author[A]{\fnms{Yi} \snm{Shen}\ead[label=e2]{ys437@cornell.edu}}
\runauthor{G. Samorodnitsky and Y. Shen}
\affiliation{Cornell University}
\address[A]{School of Operations Research\\
\quad and Information Engineering\\
Cornell University\\
220 Rhodes Hall\\
Ithaca, New York 15853\\
USA\\
\printead{e1}\\
\phantom{E-mail:\ }\printead*{e2}}
\end{aug}

\received{\smonth{10} \syear{2011}}
\revised{\smonth{6} \syear{2012}}

%
\begin{abstract}
It is, perhaps, surprising that the location of the unique supremum of
a stationary process on an interval can fail to be uniformly
distributed over that interval. We show that this distribution is
absolutely continuous in the interior of the interval and describe very
specific conditions the density has to satisfy. We establish universal
upper bounds on the density and demonstrate their optimality.
\end{abstract}

%


\end{frontmatter}
%
\section{Introduction}\label{secIntro}

The extremes of stationary processes, especially of Gaussian
processes, have attracted significant interest for a long time.
Many results are described in the books
\citet{adlertaylor2007} and \citet{azaiswschebor2009}, with
shorter versions in \citet{adler1990} and
\citet{azaiswschebor2002}. Roughly speaking, these results can be
categorized as follows: the exact distributions of the suprema
have been calculated for several particular processes; bounds on
the supremum distribution have been obtained for a large number of
processes; the asymptotic behavior of the level crossing
probability has been studied for a larger number of processes.
Almost without exception, however, these results deal with the
value of the supremum, while very little is known about the random
location of the supremum.

The present work arises from an obvious attempt to understand the
effect of stationarity of the process on the distribution of the
location of the supremum. Therefore in this paper, we look at
stationary stochastic processes in continuous, one-dimensional
time and we will consider the location of its global supremum
over a compact interval. It turns out that answering even this,
apparently simple question leads to unexpected insights.

We now discuss our setup more formally.
Let $\BX= (X(t), t\in\bbr)$ be a stationary process. If the sample
paths of the process are upper semi-continuous, then the process is
bounded from above on
any compact interval $[0,T]$, and its supremum over that interval is
attained. We are interested in the location of that supremum within
the interval $[0,T]$.

It is, of course, entirely possible that the supremum of the process
in the interval $[0,T]$ is not unique (i.e., that it is achieved at
more than one point). In that case one could be more specific and
take, for example, the left-most point in which the largest value
over the interval is achieved, as the location of the supremum. In
this paper we will sometimes deal with the situation in which, on an
event
of probability 1, the supremum is achieved at a single point. In
either case it is easy to check that the location of the supremum is a
well defined random variable.

Will the stationarity of the process guarantee a uniform
distribution of the location of the supremum over the interval?
The answer is negative. The examples in Section 9.4 of
\citet{leadbetterlindgrenrootzen1983} show that even in the case
of Gaussian processes with a uniquely attained supremum (thus
eliminating a possible bias resulting from taking the leftmost
supremum location), the supremum can still be located, with a
positive probability, at one of the endpoints of the interval and,
furthermore, the remaining mass in the interior of the interval
does not have to be uniformly distributed there.

It is, of course, the endpoints of the interval that are
responsible for the lack of uniformity. In a sense, the points
near the ends of the interval have ``fewer local competitors'' for
being the supremum than the points further from the endpoints do.
But exactly how far from having the uniform distribution can the
location of the supremum be? In this paper we give a very detailed
answer to this question by showing that this distribution is
absolutely continuous in the interior of the interval and
describing very specific conditions its density must satisfy.
This is done in Section~\ref{secassumptions}. Our results turn
out to be quite complete. In fact, we show in a companion paper,
\citet{samorodnitskyshen2012}, that for a very broad class of
stationary processes with a uniquely achieved supremum, our
description actually gives all possible distributions of its
location. In the present paper we start with treating a general
upper semi-continuous stationary process and (with one exception)
allowing the process to have multiple supremum locations within an
interval. We proceed with establishing extra conditions the
density has to satisfy if the process satisfies certain
assumptions. In Section~\ref{secupperbounds} we provide the
sharpest possible universal upper bounds on the density, both in
the general case and in the case of time-reversible stationary
processes.\looseness=1

\section{Notation and assumptions on the stationary
process} \label{secassumptions}

For the remainder of this paper $\BX= (X(t), t\in\bbr)$ is a
stationary process with upper semi-continuous sample paths, defined on
some
probability space $( \Omega, {\mathcal F}, P)$. For a compact interval
$[a,b]$, we will
denote by
\[
\tau_{\BX,[a,b]}= \min \Bigl\{ t\in[a,b]\dvtx X(t) =\sup_{a\leq
s\leq
b}X(s)
\Bigr\}.
\]
That is, $\tau_{\BX,[a,b]}$ is the first time the overall supremum in
the interval $[a,b]$ is achieved. It is elementary to check that
$\tau_{\BX,([a,b])}$ is a well-defined random variable. If $a=0$, we will
use the single variable notation $\tau_{\BX,b}$.

We denote by $F_{\BX,[a,b]}$ the law of $\tau_{\BX,[a,b]}$; it is a
probability measure on the interval $[a,b]$. If $a=0$, we have the
corresponding single variable notation $F_{\BX,b}$. The following
statements are obvious.
%
\begin{lemma} \label{llocationvar}
\textup{(i)} For any $\Delta\in\bbr$,
\[
F_{\BX,[\Delta,T+\Delta]}(\cdot) = F_{\BX,T}(\cdot-\Delta) .\vspace*{-9pt}
\]
\begin{longlist}[(ii)]
\item[(ii)] For any intervals $[c,d]\subseteq[a,b]$,
\[
F_{\BX,[a,b]}(B)\leq F_{\BX,[c,d]}(B) \qquad\mbox{for any Borel set $B
\subset[c,d]$.}
\]
\end{longlist}
\end{lemma}

The discussion of the leftmost supremum location
$\tau_{\BX,[a,b]}$ in the sequel applies equally well to the
rightmost supremum location, for instance, by considering the
time-reversed stationary process $(X(-t), t\in\bbr)$. In some
cases we will find it convenient to assume that the supremum is
achieved at a unique location. Formally, for $T>0$ we denote by
$X_*(T) = \sup_{0\leq t\leq T}X(t)$ the largest value of the
process in the interval $[0,T]$, and consider the set\looseness=1
\[
\Omega_T = \bigl\{ \omega\in\Omega\dvtx X(t_i)=X_*(T)
\mbox{ for at least two different $t_1,t_2\in[0,T]$}
\bigr\}.
\]\looseness=0
It is easy to see that $\Omega_T$ is a measurable set. The
following assumption says that, on a set of probability 1, the
supremum over interval $[0,T]$ is uniquely achieved.

\renewcommand{\theass}{$\mathrm{U}_T$}
\begin{ass}\label{assu}
$P(\Omega_T)=0$.
\end{ass}

In our previous notation, under Assumption~\ref{assu}, $\tau_{\BX,[a,b]}$ is the unique point at which the supremum
over the interval $[0,T]$ is achieved, and
$F_{\BX,T}$ is the law of that point.

Even though many of our results do not require it, the
most complete description of the distribution of the location of the
supremum that we have requires the following, additional, assumption.

\renewcommand{\theass}{$\mathrm{L}$}
\begin{ass}\label{assl}
\[
K:=\lim_{\vep\downarrow0}\frac{P ( \mbox{$\BX$ has a local
maximum in $(0,\vep)$} )}{\vep}<\infty.
\]
\end{ass}

It is easy to check that the limit in Assumption~\ref{assl} exists.
If, for example, the process $\BX$ has differentiable sample paths,
then a sufficient condition for Assumption~\ref{assl} is that the expected
number of times the process $Y(t)=X^\prime(t), t\in\bbr$ crosses
zero in a unit time interval is finite; the latter can be checked
using, for instance, Theorem 7.2.4 in
\citet{leadbetterlindgrenrootzen1983}.

Assumption~\ref{assl} rules out existence of ``too frequent'' local extrema
of the sample paths. For sample continuous processes this also rules
out rapid oscillation of the sample paths
possessed, for instance, by the Gaussian Ornstein--Uhlenbeck process of
Example~\ref{exO-U} below. In fact, we will presently see
that, at least for sample continuous processes,
under Assumption~\ref{assl} the process has, with probability 1, sample paths
of locally bounded variation.
%
\begin{lemma} \label{lKbddvar}
Let $\BX= (X(t), t\in\bbr)$ be a stationary sample upper
semi-continuous process
satisfying Assumption~\ref{assl}. Then, for any $T>0$, on an event of
probability 1 the process has finitely many local maxima and minima in
the interval $(0,T)$. In particular, if the process is sample
continuous, then its sample paths are, on event of probability 1, of
locally bounded variation.
\end{lemma}
\begin{pf}
For notational simplicity we take $T=1$. For $n=1,2,\ldots$ let
\[
N_n = \sum_{i=1}^{2^n} \one
\biggl( \mbox{a point in }  \biggl[\frac{i-1}{2^n},\frac{i}{2^n} \biggr) \mbox{ is a
local maximum of }\BX\biggr).
\]
Clearly, the sequence $N_n$ is nondecreasing, and $N_n\to N_\infty$,
where $N_\infty$ is the total number of local maxima of $\BX$ in the
interval $[0,1)$. By the monotone convergence theorem,
\begin{eqnarray*}
EN_\infty&=& \lim_{n\to\infty} EN_n
\\
&\leq&\limsup_{n\to\infty} 2^n P \bigl( \mbox{$\BX$ has a local
maximum in $\bigl(0,2^{-n}\bigr)$} \bigr) \leq K.
\end{eqnarray*}
Therefore, $N_\infty<\infty$ a.s. Since between any two distinct local
minima there is a local maximum, the number of local minima in $[0,1)$
is a.s. finite as well. Since a sample continuous process must have
a monotone path between any two consecutive local extrema, the
lemma has been proved.
\end{pf}

\section{Description of the possible distributions of the location of
the supremum} \label{secdensitydescr}

We start with a result showing existence of a density in the interior
of the interval $[0,T]$ of the leftmost location\vadjust{\goodbreak} of the supremum in that
interval for any upper semi-continuous stationary process, as well as
conditions this density has to satisfy. Only one of the statements of
the theorem requires Assumption~\ref{assu}, in which case the statement
applies to the unique location of the supremum. See Remark
\ref{rkzerodensity} in the sequel.
%
\begin{theorem} \label{tdensityrcll}
Let $\BX= (X(t), t\in\bbr)$ be a stationary sample upper
semi-continuous process. Then the restriction of the law
$F_{\BX,T}$ to the interior $(0,T)$ of the interval is absolutely
continuous. The density, denoted by $f_{\BX,T}$, can be taken to be
equal to the right derivative of the cdf $F_{\BX,T}$, which exists at
every point in the interval $(0,T)$. In this case the density is right
continuous, has left limits, and has the following properties:
\begin{longlist}[(a)]
\item[(a)] The limits
\[
f_{\BX,T}(0+)=\lim_{t\to0} f_{\BX,T}(t) \quad\mbox{and}\quad
f_{\BX
,T}(T-)=\lim_{t\to T} f_{\BX,T}(t)
\]
exist.

\item[(b)] The density has a universal upper bound given by
%
\begin{equation}
\label{edensitybound} f_{\BX,T}(t) \leq\max \biggl(\frac1t,
\frac{1}{T-t} \biggr),\qquad 0<t<T .
\end{equation}

\item[(c)] Assume that the process satisfies Assumption~\ref{assu}. Then
the density is bounded away from zero,
%
\begin{equation}
\label{edensitylowerbound} \inf_{0<t<T}f_{\BX,T}(t) >0 .
\end{equation}

\item[(d)] The density has a bounded variation away from the endpoints of
the interval. Furthermore, for every $0<t_1<t_2<T$,
%
\begin{equation}
\label{eTVaway} \qquad\mathrm{TV}_{(t_1, t_2)}(f_{\BX,T}) \leq\min \bigl(
f_{\BX,T}(t_1), f_{\BX,T}(t_1-) \bigr) +
\min \bigl( f_{\BX,T}(t_2), f_{\BX,T}(t_2-)
\bigr) ,
\end{equation}
where
\[
\mathrm{TV}_{(t_1, t_2)}(f_{\BX,T}) = \sup\sum_{i=1}^{n-1}
\bigl|f_{\BX,T}(s_{i+1})- f_{\BX,T}(s_i) \bigr|
\]
is the total variation of $f_{\BX,T}$ on the interval $(t_1,t_2)$, and
the supremum is taken over all choices of $t_1<s_1<\cdots<s_n<t_2$.

\item[(e)] The density has a bounded positive variation at the left
endpoint and a bounded negative variation at the right
endpoint. Furthermore, for every $0<\vep<T$,
%
\begin{equation}
\label{eTVleft} \mathrm{TV}^+_{(0,\vep)}(f_{\BX,T}) \leq\min \bigl(
f_{\BX,T}(\vep), f_{\BX,T}(\vep-) \bigr)
\end{equation}
and
%
\begin{equation}
\label{eTVright} \mathrm{TV}^-_{(T-\vep,T)}(f_{\BX,T}) \leq\min \bigl(
f_{\BX,T}(T-\vep), f_{\BX,T}(T-\vep-) \bigr) ,
\end{equation}
where for any interval $0\leq a<b\leq T$,
\[
\mathrm{TV}^\pm_{(a,b)}(f_{\BX,T}) = \sup\sum
_{i=1}^{n-1} \bigl(f_{\BX,T}(s_{i+1})-
f_{\BX,T}(s_i) \bigr)_\pm
\]
is the positive (negative) variation of $f_{\BX,T}$ on the interval
$(a,b)$, and the supremum is taken over all choices of
$a<s_1<\cdots<s_n<b$.

\item[(f)] The limit $f_{\BX,T}(0+)<\infty$ if
and only if $\mathrm{TV}_{(0,\vep)}(f_{\BX,T}) <\infty$ for some (equivalently,
any) $0<\vep<T$, in which case
%
\begin{equation}
\label{eTVleft1} \mathrm{TV}_{(0,\vep)}(f_{\BX,T}) \leq f_{\BX,T}(0+)+
\min \bigl( f_{\BX
,T}(\vep), f_{\BX,T}(\vep-) \bigr) .
\end{equation}
Similarly, $f_{\BX,T}(T-)<\infty$ if and only if
$\mathrm{TV}_{(T-\vep,T)}(f_{\BX,T}) <\infty$ for some (equivalently,
any) $0<\vep<T$, in which case
%
\begin{equation}
\label{eTVright1} \mathrm{TV}_{(T-\vep,T)}(f_{\BX,T}) \leq\min \bigl(
f_{\BX,T}(T-\vep), f_{\BX,T}(T-\vep-) \bigr) + f_{\BX,T}(T-)
.
\end{equation}
\end{longlist}
\end{theorem}
\begin{pf}
Choose $0<\delta<T/2$. We claim that for every $\delta\leq t\leq
T-\delta$, for every $\rho>0$ and every $0<\vep<\delta\rho/(1+\rho)$,
%
\begin{equation}
\label{eabsctybound} P ( t<\tau_{\BX,T}\leq t+\vep ) \leq \vep(1+\rho)
\max \biggl(\frac1t,\frac{1}{T-t} \biggr) .
\end{equation}
This statement, once proved, will imply absolute continuity of
$F_{\BX,T}$ on the interval $(\delta, T-\delta)$ and, since $\delta>0$
can be taken to be arbitrarily small, also on $(0,T)$. Further,
\eqref{eabsctybound} will imply that the version of the density
given by
\[
f_{\BX,T}(t) = \limsup_{\vep\downarrow0} \frac1\vep P ( t<
\tau_{\BX,T}\leq t+\vep ),\qquad 0<t<T ,
\]
satisfies bound \eqref{edensitybound}.

We proceed to prove \eqref{eabsctybound}. Suppose that, to the
contrary, \eqref{eabsctybound} fails for some $\delta\leq t\leq
T-\delta$ and $0<\vep<\delta\rho/(1+\rho)$. Choose
\[
\vep<\theta<\frac{\rho}{1+\rho}\delta
\]
and $0<a<t<b<T$ such that
\[
\min ( t,T-t ) - \theta<b-a<\min ( t,T-t )-\vep.
\]
For $a\leq s\leq b$, by stationarity, we have
%
\begin{equation}
\label{eeachs} P ( s<\tau_{\BX,[s-t,s-t+T]}\leq s+\vep ) > \vep(1+\rho)\max \biggl(
\frac1t,\frac{1}{T-t} \biggr) .
\end{equation}
Further, let $a\leq s_1<s_1+\vep\leq s_2\leq b$. We check next that
%
\begin{equation}
\label{edisjoint} \{ s_j<\tau_{\BX,[s_j-t,s_j-t+T]}\leq s_j+
\vep, j=1,2 \}=\varnothing.
\end{equation}
Indeed, let $\Omega_{s_1,s_2}$ be the event in
\eqref{edisjoint}. Note that the intervals $(s_1,s_1+\vep)$ and
$(s_2,s_2+\vep)$ are disjoint and, by the choice of the parameters $a$
and $b$, each of these two intervals is a subinterval of both
$[s_1-t,s_1-t+T]$ and $[s_2-t,s_2-t+T]$. Therefore, on the event
$\Omega_{s_1,s_2}$ we cannot have
\[
X (\tau_{\BX,[s_1-t,s_1-t+T]} )< X (\tau_{\BX,[s_2-t,s_2-t+T]} )
\]
for otherwise $\tau_{\BX,[s_1-t,s_1-t+T]}$ would fail to be a
location of the maximum over the interval $[s_1-t,s_1-t+T]$. For the
same reason, on the event
$\Omega_{s_1,s_2}$ we cannot have
\[
X (\tau_{\BX,[s_1-t,s_1-t+T]} )> X (\tau_{\BX,[s_2-t,s_2-t+T]} ) .
\]
Finally, on the event
$\Omega_{s_1,s_2}$ we cannot have
\[
X (\tau_{\BX,[s_1-t,s_1-t+T]} )= X (\tau_{\BX,[s_2-t,s_2-t+T]} )
\]
for otherwise $\tau_{\BX,[s_2-t,s_2-t+T]}$ would fail to be the
leftmost location of the maximum over the interval $[s_2-t,s_2-t+T]$.
This establishes \eqref{edisjoint}.

We now apply \eqref{eeachs} and \eqref{edisjoint} to the points
$s_i = a+i\vep, i=0,1,\ldots, \lceil(b-a)/\vep\rceil-1$. We have
\begin{eqnarray*}
1&\geq& P \Biggl( \bigcup_{i=0}^{\lceil(b-a)/\vep\rceil-1} \{
s_i<\tau_{\BX,[s_i-t,s_i-t+T]}\leq s_i+\vep \} \Biggr)
\\
&= &\sum_{i=0}^{\lceil(b-a)/\vep\rceil-1} P ( s_i<
\tau_{\BX,[s_i-t,s_i-t+T]}\leq s_i+\vep ) \\
&>& \frac{b-a}{\vep} \vep(1+\rho)
\max \biggl(\frac1t,\frac
{1}{T-t} \biggr)
\\
&>& \bigl( \min ( t,T-t ) - \theta \bigr) (1+\rho)\max \biggl(\frac1t,
\frac{1}{T-t} \biggr)
\\
&>& \biggl( 1- \frac{\delta}{\min( t,T-t)} \frac{\rho}{1+\rho
} \biggr) (1+\rho) \geq \biggl(
1- \frac{\rho}{1+\rho} \biggr) (1+\rho) =1
\end{eqnarray*}
by the choice of $\theta$. This contradiction proves
\eqref{eabsctybound}.

Before proceeding with the proof of Theorem~\ref{tdensityrcll}, we
pause to prove the following important lemma.
%
\begin{lemma} \label{lkey}
Let $0\leq\Delta<T$. Then for every $0\leq\delta\leq\Delta$,
$f_{\BX,T-\Delta}(t)\geq f_{\BX,T}(t+\delta)$ almost everywhere in
$(0,T-\Delta)$. Furthermore, for every such $\delta$ and every
$\vep_1,\vep_2\geq0$, such that $\vep_1+\vep_2<T-\Delta$,
%
\begin{eqnarray}
\label{edensityvarint} &&\int_{\vep_1}^{T-\Delta-\vep_2} \bigl(
f_{\BX,T-\Delta}(t)-f_{\BX,T}(t+\delta) \bigr) \,dt
\nonumber
\\[-8pt]
\\[-8pt]
\nonumber
&&\qquad\leq \int_{\vep_1}^{\vep_1+\delta} f_{\BX,T}(t) \,dt +
\int_{T-\Delta-\vep_2+\delta}^{T-\vep_2} f_{\BX,T}(t) \,dt .
\end{eqnarray}
\end{lemma}
\begin{pf}
We simply use Lemma~\ref{llocationvar}. For any Borel set
$B\subseteq(0,T-\Delta)$ we have
\begin{eqnarray*}
\int_B f_{\BX,T-\Delta}(t) \,dt &=& P ( \tau_{\BX,T-\Delta}
\in B ) \geq P ( \tau_{\BX,[-\delta,T-\delta]}\in B )
\\
&= &\int_B f_{\BX,[-\delta,T-\delta]}(t) \,dt = \int
_B f_{\BX,T}(t+\delta) \,dt ,
\end{eqnarray*}
which shows that $f_{\BX,T-\Delta}(t)\geq f_{\BX,T}(t+\delta)$ almost
everywhere in $(0,T-\Delta)$.

For \eqref{edensityvarint}, notice that by Lemma
\ref{llocationvar},
\begin{eqnarray*}
&&\int_{\vep_1}^{T-\Delta-\vep_2} \bigl( f_{\BX,T-\Delta}(t)-f_{\BX,T}(t+
\delta) \bigr) \,dt
\\
&&\qquad= P \bigl( \tau_{\BX,T-\Delta}\in(\vep_1,T-\Delta-
\vep_2) \bigr) - P \bigl( \tau_{\BX,T}\in(\vep_1+
\delta,T-\Delta-\vep_2+\delta ) \bigr)
\\
&&\qquad= P \bigl( \tau_{\BX,T}\notin (\vep_1+\delta,T-\Delta-
\vep_2+\delta) \bigr) - P \bigl( \tau_{\BX,T-\Delta}\notin(
\vep_1,T-\Delta-\vep_2) \bigr)
\\
&&\qquad= P \bigl( \tau_{\BX,T}\in[0,\vep_1+\delta) \bigr) + P
\bigl( \tau_{\BX,T}\in(T-\Delta-\vep_2+\delta, T] \bigr)
\\
&&\qquad\quad{}- P \bigl( \tau_{\BX,T-\Delta}\in[0,\vep_1) \bigr) - P (
\tau_{\BX,T-\Delta}\in(T-\Delta-\vep_2, T-\Delta ] \bigr)
\\
&&\qquad= P \bigl( \tau_{\BX,T}\in(\vep_1,\vep_1+
\delta) \bigr) + \bigl( P \bigl( \tau_{\BX,T}\in[0,\vep_1
) \bigr)- P \bigl( \tau_{\BX,T-\Delta}\in[0,\vep_1) \bigr) \bigr)
\\
&&\qquad\quad{}+ P \bigl( \tau_{\BX,T}\in(T-\Delta-\vep_2+\delta, T-
\vep_2) \bigr)
\\
&&\qquad\quad{}+ \bigl( P\bigl ( \tau_{\BX,T}\in(T-\vep_2,T] \bigr) - P \bigl(
\tau_{\BX,[\Delta,T]}\in(T-\vep_2,T] \bigr) \bigr)
\\
&&\qquad\leq P \bigl( \tau_{\BX,T}\in(\vep_1,\vep_1+
\delta) \bigr) + P \bigl( \tau_{\BX,T}\in(T-\Delta-\vep_2+
\delta, T-\vep_2) \bigr)
\\
&&\qquad= \int_{\vep_1}^{\vep_1+\delta} f_{\BX,T}(t) \,dt + \int
_{T-\Delta-\vep_2+\delta}^{T-\vep_2} f_{\BX,T}(t) \,dt
\end{eqnarray*}
as required.
\end{pf}
We return now to the proof of Theorem~\ref{tdensityrcll}. Our next
goal is to prove that the cdf $F_{\BX,T}$ is right differentiable at
every point in the interval $(0,T)$. Since we already know that
$F_{\BX,T}$ is absolutely continuous on $(0,T)$, the set
%
\begin{equation}
\label{eA} A = \bigl\{ t\in(0,T)\dvtx F_{\BX,T} \mbox{ is not right
differentiable at $t$} \bigr\}
\end{equation}
has Lebesgue measure zero. Define next
%
\begin{equation}
\label{eB}   B= \bigl\{ t\in A^c\dvtx f_{\BX,T} \mbox{ restricted
to $A^c$ does not have a right limit at $t$} \bigr\}.\hspace*{-35pt}
\end{equation}
We claim that the set $B$ is at most countable. To see this, we define
for $ t\in A^c$
\[
L(t) = \limsup_{s\downarrow t, s\in A^c} f_{\BX,T}(s),\qquad l(t) = \liminf_{s\downarrow t, s\in A^c}
f_{\BX,T}(s) .
\]
Our claim about set $B$ will follow once we check that for any
$0<\vep<T/2$ and $\theta>0$, the set
\[
B_{\vep,\theta}= \bigl\{ t\in A^c\cap(\vep,T-\vep)\dvtx L(t)-l(t)>
\theta \bigr\}
\]
is finite. In fact, we will show that the cardinality of
$B_{\vep,\theta}$ cannot be larger than $4/(\vep\theta)$. If not, let
$N>4/(\vep\theta)$ and find points $\vep<t_1<t_2<\cdots
<t_N<T-\vep$. Choose $\delta>0$ so small that $\delta<\vep/2$ and
\[
0<\delta<\tfrac12 \min ( t_1-\vep,t_2-t_1,
\ldots, t_N-t_{N_1}, T-\vep-t_N ) .
\]
Let now $i=1,\ldots, N$ and choose a sequence $s_n\downarrow t_i,
s_n\in A^c$, such that\break $f_{\BX,T}(s_n)\to L(t_i)$. Consider $n$ so
large that $s_n-t_i<\delta/3$, and let
\[
j\geq\frac{3}{\delta-(s_n-t_i)}
\]
be an integer. We have
\[
P \bigl( \tau_{\BX,T-\delta}\in(t_i-\delta,t_i) \bigr)
\geq\sum_{k=0}^{\lfloor j(\delta-(s_n-t_i))\rfloor-1} P \bigl(
\tau_{\BX,T-\delta}\in\bigl(t_i-(k+1)/j,t_i-k/j\bigr)
\bigr) ,
\]
and for each $k$ as in the sum
\[
h_k:= s_n-t_i+\frac{k+1}{j}\in(0,
\delta] .
\]
Therefore, by Lemma~\ref{llocationvar}
\begin{eqnarray*}
&&P \bigl( \tau_{\BX,T-\delta}\in(t_i-\delta,t_i) \bigr)
\\
&&\qquad\geq\sum_{k=0}^{\lfloor j(\delta-(s_n-t_i))\rfloor-1} P \bigl(
\tau_{\BX,T}\in\bigl(t_i-(k+1)/j+h_k,t_i-k/j+h_k
\bigr) \bigr)
\\
&&\qquad= \bigl\lfloor j\bigl(\delta-(s_n-t_i)\bigr)\bigr
\rfloor P \bigl( \tau_{\BX,T}\in (s_n,s_n+1/j)
\bigr) \to\bigl(\delta-(s_n-t_i)\bigr)
f_{\BX,T}(s_n)
\end{eqnarray*}
as $j\to\infty$. Letting $n\to\infty$, we conclude that
%
\begin{equation}
\label{eproblarge} P \bigl( \tau_{\BX,T-\delta}\in(t_i-
\delta,t_i) \bigr) \geq\delta L(t_i),\qquad i=1,\ldots, N .
\end{equation}

Similarly, for $i=1,\ldots, N$ choose a sequence $w_n\downarrow t_i,
w_n\in A^c$, such that $f_{\BX,T}(w_n)\to l(t_i)$. For large $n$ and
$j$ we have
\begin{eqnarray*}
&&P \bigl( \tau_{\BX,T+\delta}\in(t_i,t_i+\delta) \bigr)\\
&&\qquad= P \bigl( \tau_{\BX,T+\delta}\in(t_i,w_n) \bigr) + P
\bigl( \tau_{\BX,T+\delta}\in(w_n,w_n+\delta) \bigr)
\\
&&\qquad\leq P \bigl( \tau_{\BX,T+\delta}\in(t_i,w_n) \bigr)+
\sum_{k=0}^{\lceil\delta j\rceil-1} P \bigl( \tau_{\BX,T+\delta}
\in\bigl(w_n+k/j,w_n+(k+1)/j\bigr) \bigr) .
\end{eqnarray*}
For each $k$ as in the sum above,
\[
h_k:=\frac{k}{j}\in[0,\delta] .
\]
Therefore, by Lemma~\ref{llocationvar},
\begin{eqnarray*}
&&P \bigl( \tau_{\BX,T+\delta}\in(t_i,t_i+\delta) \bigr)
\\
&&\qquad \leq P \bigl( \tau_{\BX,T+\delta}\in(t_i,w_n) \bigr)+
\lceil\delta j\rceil P \bigl( \tau_{\BX,T}\in(w_n,w_n+1/j)
\bigr).
\end{eqnarray*}
Letting, once again, first $j\to\infty$ and then $n\to\infty$, we
conclude that
%
\begin{equation}
\label{eprobsmall} P \bigl( \tau_{\BX,T+\delta}\in(t_i,t_i+
\delta) \bigr) \leq\delta l(t_i),\qquad i=1,\ldots, N .
\end{equation}
Now we use the estimate in Lemma~\ref{lkey} as follows. By the
definition of the point~$t_i$ and the smallness of $\delta$,
\begin{eqnarray*}
N\delta\theta&\leq& P \Biggl( \tau_{\BX,T-\delta}\in\bigcup
_{i=1}^N (t_i-\delta,t_i)
\Biggr) - P \Biggl( \tau_{\BX,T+\delta}\in\bigcup_{i=1}^N
(t_i,t_i+\delta) \Biggr)
\\
&= &\int_{\bigcup_{i=1}^N (t_i-\delta,t_i)} \bigl( f_{\BX,T-\delta
}(t) - f_{\BX,T+\delta}(t+
\delta) \bigr) .
\end{eqnarray*}
Using the fact that
\[
\bigcup_{i=1}^N (t_i-
\delta,t_i) \subset(\vep-\delta, T-\vep) ,
\]
and that, by Lemma~\ref{lkey}, the integrand above is
a.e. nonnegative, we have by the estimate in that lemma that the
integral above does not exceed
\begin{eqnarray*}
&&\int_{\vep-\delta}^{T-\vep} \bigl( f_{\BX,T-\delta}(t) -
f_{\BX,T+\delta}(t+\delta) \bigr) \,dt
\\
&&\qquad\leq\int_{\vep-\delta}^\vep f_{\BX,T+\delta}(t) \,dt + \int
_{T-\vep+\delta}^{T-\vep+2\delta} f_{\BX,T+\delta}(t) \,dt .
\end{eqnarray*}
Applying the already proved \eqref{edensitybound}, we conclude that
\[
N\delta\theta\leq2\frac{\delta}{\vep-\delta}\leq\frac{4\delta
}{\vep} ,
\]
and this contradicts the assumption that we can choose
$N>4/(\vep\theta)$. This proves that the set $B$ in \eqref{eB} is at
most countable. We notice, further, that
%
\begin{eqnarray}
\label{enotAB} f_{\BX,T}(t) &=& \lim_{s\downarrow t} \frac{1}{s-t} P (
t<\tau_{\BX,T}\leq s )
\nonumber
\\[-8pt]
\\[-8pt]
\nonumber
&=& \lim_{s\downarrow t} \frac{1}{s-t} \int_t^s
f_{\BX,T}(w) \,dw = \lim_{w\downarrow t, w\in A^c\setminus B}f_{\BX,T}(w)
\end{eqnarray}
for every $t\in A^c\setminus B$ [recall the set $A$ is defined in
\eqref{eA}].\vadjust{\goodbreak}

Now we are ready to prove that the right derivative of the cdf
$F_{\BX,T}$ exists at every point in the interval $(0,T)$.
Suppose, to the contrary, that this is not so. Then there is $t\in
(0,T)$ and real numbers $a<b$ such that
\[
\liminf_{\vep\downarrow0}\frac{F_{\BX,T}(t+\vep)-F_{\BX
,T}(t)}{\vep} <a<b<\limsup_{\vep\downarrow
0}
\frac{F_{\BX,T}(t+\vep)-F_{\BX,T}(t)}{\vep}.
\]
This implies that there is a sequence $t_n\downarrow t$ with $t_n\in
A^c\setminus B$ for each $n$ such that
\[
f_{\BX,T}(t_{2n-1})>b, \qquad f_{\BX,T}(t_{2n})<a\qquad
\mbox{for all $n=1,2,\ldots.$}
\]
We can and will choose $t_1$ so close to $t$ that $t_1<(T+t)/2$.

Notice that by \eqref{enotAB}, for every $n=1,2,\ldots$ there is
$\delta_n>0$ such that
\begin{eqnarray*}
f_{\BX,T}(w)&>&b \qquad \mbox{a.e. in $(t_{2n-1}, t_{2n-1}+
\delta_{2n-1})$} ,
\\
f_{\BX,T}(w)&<& a\qquad  \mbox{a.e. in $(t_{2n}, t_{2n}+
\delta_{2n})$}
\end{eqnarray*}
for $n=1,2,\ldots.$

Let now $m\geq1$, and consider $s>0$ so small that both
$s<\min_{n=1,\ldots, 2m}\delta_n$ and $t_1<(T+t)/2-s$. Observe that
\begin{eqnarray*}
&&\int_t^{(T+t)/2} \bigl( f_{\BX,T}(w+s) -
f_{\BX,T}(w) \bigr)_+ \,dw
\\
&&\qquad\geq\int_t^{t+s} \sum
_{i=0}^{\lfloor(T-t)/2s\rfloor-1} \bigl( f_{\BX,T}\bigl(w+(i+1)s
\bigr) - f_{\BX,T}(w+is) \bigr)_+ \,dw ,
\end{eqnarray*}
and for every point $w\in(t,t+s)$, each one of the intervals
$(t_n,t_n+\delta_n), n=1,\ldots, 2m$, contains at least one of the
points in the
finite sequence $w+is, i =0,1,\ldots, \lfloor(T-t)/2s\rfloor
-1$. By construction, apart from a set of points $w\in(t,t+s)$ of
measure zero, those points of the kind $w+is$ that fall in the
odd-numbered intervals satisfy $f_{\BX,T}(w+is)>b$, and those points
that fall in the
even-numbered intervals satisfy $f_{\BX,T}(w+is)<a$. We conclude that
\[
\sum_{i=0}^{\lfloor(T-t)/2s\rfloor-1} \bigl( f_{\BX,T}
\bigl(w+(i+1)s\bigr) - f_{\BX,T}(w+is) \bigr)_+\geq m(b-a)
\]
a.e. in $(t,t+s)$. Therefore, for all $s>0$ small enough,
\[
\int_t^{(T+t)/2} \bigl( f_{\BX,T}(w+s) -
f_{\BX,T}(w) \bigr)_+ \,dw \geq s m(b-a)
\]
and, since $m$ can be taken arbitrarily large, we conclude that
%
\begin{equation}
\label{eaveragelarge} \lim_{s\downarrow0} \frac1s \int_t^{(T+t)/2}
\bigl( f_{\BX
,T}(w+s) - f_{\BX,T}(w) \bigr)_+ \,dw = \infty.
\end{equation}
We will see that this is, however, impossible, and the resulting
contradiction will prove that the right derivative of the cdf
$F_{\BX,T}$ exists at every point in the interval $(0,T)$.

Indeed, recall that by Lemma~\ref{lkey}, for all $s>0$ small enough,
\[
f_{\BX,T-2s}(w-s)\geq f_{\BX,T}(w+s) \qquad\mbox{a.e. on $(s,T-s)\supset
\bigl(t, (T+t)/2\bigr).$}
\]
Therefore, for such $s$,
\begin{eqnarray*}
&&\int_t^{(T+t)/2} \bigl( f_{\BX,T}(w+s) -
f_{\BX,T}(w) \bigr)_+ \,dw
\\
&&\qquad\leq \int_t^{(T+t)/2} \bigl( f_{\BX,T-2s}(w-s) -
f_{\BX,T}(w) \bigr)_+ \,dw
\\
&&\qquad\leq\int_{t-s}^{(T+t)/2-s} \bigl( f_{\BX,T-2s}(w) -
f_{\BX,T}(w+s) \bigr) \,dw
\end{eqnarray*}
since, by another application of Lemma~\ref{lkey}, the integrand is
a.e. nonnegative over the range of integration. Applying
\eqref{edensityvarint}, we see that
\begin{eqnarray*}
&&\int_t^{(T+t)/2} \bigl( f_{\BX,T}(w+s) -
f_{\BX,T}(w) \bigr)_+ \,dw
\\
&&\qquad\leq\int_{t-s}^t f_{\BX,T}(w) \,dw + \int
_{(T+t)/2}^{(T+t)/2+s} f_{\BX,T}(w) \,dw .
\end{eqnarray*}
However, we already know that the density $ f_{\BX,T}$ is bounded on
any subinterval of $(0,T)$ that is bounded away from both
endpoints. Therefore, the upper bound obtained above shows that
\eqref{eaveragelarge} is impossible. Hence the existence of the
right derivative everywhere, which then coincides with the version of
the density $f_{\BX,T}$ chosen above.

Next we check that this version of the density is right continuous. To
this end we recall that we already know that the set $A$ in
\eqref{eA} is empty. Next, we
rule out existence of a point $t\in(0,T)$ such the limit of
$f_{\BX,T}(s)$ as $s\downarrow t$ over $s\in B^c$ does not exist.
Suppose that, to the
contrary, that such $t$ exists. This means that there are real
numbers $a<b$ and a sequence $t_n\downarrow t$ with $t_n\in
B^c$ for each $n$ such that
\[
f_{\BX,T}(t_{2n-1})>b, \qquad f_{\BX,T}(t_{2n})<a\qquad
\mbox{for all $n=1,2,\ldots.$}
\]
However, we have already established that such a sequence cannot
exist.

As in \eqref{enotAB}, we see that for every $t\in(0,T)$
\[
f_{\BX,T}(t) = \lim_{s\downarrow t, s\in B^c}f_{\BX,T}(s),
\]
and since the set $B$ is at most countable, the restriction to $s\in
B^c$ in the above limit statement can be removed. This proves right
continuity of the version of the density given by the right
derivative of $F_{\BX,T}$. The proof of existence of left limits is
similar.\vadjust{\goodbreak}

Next, we address the variation of the version of the density we are
working with away from the endpoints of the interval $(0,T)$.
Let $0<t_1<t_2<T$. We start with a preliminary calculation.
Let $0<r_n<T-t_2$. Introduce the notation
\begin{eqnarray*}
C_+&=& \bigl\{ t\in(t_1,t_2)\dvtx f_{\BX,T}(t+r_n)
\geq f_{\BX,T}(t) \bigr\} ,
\\
C_-&=& \bigl\{ t\in(t_1,t_2)\dvtx f_{\BX,T}(t+r_n)<
f_{\BX,T}(t) \bigr\} ,
\end{eqnarray*}
so that
\begin{eqnarray*}
&&\int_{t_1}^{t_2} \bigl| f_{\BX,T}(t+r_n)-f_{\BX,T}(t)
\bigr| \,dt
\\
&&\qquad= \int_{C_+} \bigl( f_{\BX,T}(t+r_n)-f_{\BX,T}(t)
\bigr) \,dt + \int_{C_-} \bigl( f_{\BX,T}(t)-f_{\BX,T}(t+r_n)
\bigr) \,dt .
\end{eqnarray*}
To estimate the two terms we will once again use Lemma
\ref{lkey}. Since
\[
f_{\BX,T-r_n}(t)\geq f_{\BX,T}(r_n+t) \qquad\mbox{a.e. on }
(0,T-r_n) \supset(t_1,t_2)
\]
for $n$ large enough, for such $n$, we have the upper bound
\begin{eqnarray*}
\int_{C_+} \bigl( f_{\BX,T}(t+r_n)-f_{\BX,T}(t)
\bigr) \,dt &\leq&\int_{C_+} \bigl( f_{\BX,T-r_n}(t)-f_{\BX,T}(t)
\bigr) \,dt
\\
&\leq&\int_{t_1}^{t_2} \bigl( f_{\BX,T-r_n}(t)-f_{\BX,T}(t)
\bigr) \, dt .
\end{eqnarray*}
We now once again use \eqref{edensityvarint} to conclude that for
all $n$ large, we have
\[
\int_{C_+} \bigl( f_{\BX,T}(t+r_n)-f_{\BX,T}(t)
\bigr) \,dt \leq\int_{t_2}^{t_2+r_n} f_{\BX,T}(t)
\,dt
\]
so that
\[
\limsup_{n\to\infty} \frac{1}{r_n}\int_{C_+} \bigl(
f_{\BX,T}(t+r_n)-f_{\BX,T}(t) \bigr) \,dt \leq
f_{\BX,T}(t_2) .
\]
Similarly, by Lemma~\ref{lkey},
\[
f_{\BX,T}(t+r_n)\geq f_{\BX,T+r_n}(t+r_n)
\qquad\mbox{a.e. on } (0,T-r_n) \supset(t_1,t_2)
\]
for $n$ large enough, and we obtain, for such $n$, using \eqref{edensityvarint}
\begin{eqnarray*}
\int_{C_-} \bigl( f_{\BX,T}(t)-f_{\BX,T}(t+r_n)
\bigr) \,dt& \leq&\int_{C_-} \bigl( f_{\BX,T}(t)-f_{\BX,T+r_n}(t+r_n)
\bigr) \,dt
\\
&\leq&\int_{t_1}^{t_2} \bigl( f_{\BX,T}(t)-f_{\BX,T+r_n}(t+r_n)
\bigr) \,dt\\
& \leq&\int_{t_1}^{t_1+r_n} f_{\BX,T+r_n}(t)
\,dt .
\end{eqnarray*}
This can, in turn, be bounded from above both by
\[
\int_{t_1}^{t_1+r_n} f_{\BX,T}(t) \,dt
\]
and by
\[
\int_{t_1}^{t_1+r_n} f_{\BX,T}(t-r_n)
\,dt = \int_{t_1-r_n}^{t_1} f_{\BX
,T}(t) \,dt .
\]
Therefore,
\[
\limsup_{n\to\infty} \frac{1}{r_n} \int_{C_-} \bigl(
f_{\BX,T}(t)-f_{\BX,T}(t+r_n) \bigr) \,dt \leq\min
\bigl( f_{\BX,T}(t_1), f_{\BX,T}(t_1-)
\bigr) .
\]
Overall, we have proved that
%
\begin{eqnarray}
\label{evarboundint} &&\limsup_{n\to\infty} \frac{1}{r_n} \int
_{t_1}^{t_2} \bigl| f_{\BX,T}(t+r_n)-f_{\BX,T}(t)
\bigr| \,dt
\nonumber
\\[-8pt]
\\[-8pt]
\nonumber
&&\qquad\leq\min \bigl( f_{\BX,T}(t_1), f_{\BX,T}(t_1-)
\bigr) + f_{\BX,T}(t_2) .
\end{eqnarray}

To relate \eqref{evarboundint} to the total variation of the
density $f_{\BX,T}$ over the interval $(t_1,t_2)$, we notice first
that by the right continuity of the density, it is
enough to consider the regularly spaced points $s_i=t_1+ir_n,
i=1,\ldots
, n$, where $r_n = (t_2-t_1)/(n+1)$ for some
$n=1,2,\ldots.$ Write
\[
\int_{t_1}^{t_2} \bigl| f_{\BX,T}(t+r_n)-f_{\BX,T}(t)
\bigr| \,dt = \int_{t_1}^{t_1+r_n} \sum
_{i=0}^n \bigl| f_{\BX,T}\bigl(t+(i+1)r_n
\bigr)-f_{\BX,T}(t+ir_n) \bigr| \,dt
\]
and observe that
\[
\lim_{n\to\infty}\sum_{i=0}^n \bigl|
f_{\BX,T}\bigl(t+(i+1)r_n\bigr)-f_{\BX,T}(t+ir_n)
\bigr| \geq \mathrm{TV}_{(t_1,
t_2)}(f_{\BX,T})
\]
uniformly in $t\in(t_1,t_2)$. Therefore, by \eqref{evarboundint}
\begin{eqnarray*}
&&\min \bigl( f_{\BX,T}(t_1), f_{\BX,T}(t_1-)
\bigr) + f_{\BX,T}(t_2)\\
&&\qquad\geq \limsup_{n\to\infty}
\frac{1}{r_n} \int_{t_1}^{t_2}\bigl |
f_{\BX,T}(t+r_n)-f_{\BX,T}(t) \bigr| \,dt
\\
&&\qquad\geq\limsup_{n\to\infty} \frac{1}{r_n} \int_{t_1}^{t_1+r_n}
\sum_{i=0}^n \bigl| f_{\BX,T}
\bigl(t+(i+1)r_n\bigr)-f_{\BX,T}(t+ir_n) \bigr| \,dt\\
&&\qquad\geq \mathrm{TV}_{(t_1,
t_2)}(f_{\BX,T}) .
\end{eqnarray*}
Now bound \eqref{eTVaway} follows from the obvious fact that
\[
\mathrm{TV}_{(t_1, t_2)}(f_{\BX,T}) = \lim_{\vep\downarrow0} \mathrm{TV}_{(t_1,
t_2-\vep)}(f_{\BX,T}).
\]
Furthermore, the proof of \eqref{eTVleft} and \eqref{eTVright} is
the same as the proof of \eqref{eTVaway}, with each one using one
side of the two-sided calculation performed above for
\eqref{eTVaway}.

Next, the boundedness of the positive variation of the density at
zero, clearly, implies that the limit $f_{\BX,T}(0+)=
\lim_{t\downarrow0} f_{\BX,T}(t)$ exists, while the boundedness of
the negative variation of the density at $T$ implies that the limit
$f_{\BX,T}(T-)= \lim_{t\uparrow T} f_{\BX,T}(t)$ exists as well. If
$\mathrm{TV}_{(0,\vep)}(f_{\BX,T}) <\infty$ for some $0<\vep<T$, then,
trivially, $f_{\BX,T}(0+)<\infty$. On the other hand, if
$f_{\BX,T}(0+)<\infty$, then the same argument as we used in proving
\eqref{eTVaway}, shows that for any $0<\vep<T$,
\[
\mathrm{TV}^-_{(0,\vep)}(f_{\BX,T}) \leq f_{\BX,T}(0+) ,
\]
which, together with \eqref{eTVleft}, both shows that
$\mathrm{TV}_{(0,\vep)}(f_{\BX,T}) <\infty$ and proves \eqref{eTVleft1}. One
can prove the statement of part (f) of the theorem concerning the
behavior of the density at the right endpoint of the interval in the
same way.

It only remains to prove part (c) of the theorem, namely the fact that
the version of the density given by the right derivative of the cdf
$F_{\BX,T}$ is bounded away from zero. Recall that Assumption~\ref{assu} is
in effect here.

Suppose, to the contrary, that
\eqref{edensitylowerbound} fails and introduce the notation
\begin{eqnarray*}
t_1&=& \inf \Bigl\{ s\in(0,T)\dvtx\inf_{0<t<s}
f_{\BX,T}(t)=0 \Bigr\},
\\
t_2&=& \sup \Bigl\{ s\in(0,T)\dvtx\inf_{s<t<T}
f_{\BX,T}(t)=0 \Bigr\} .
\end{eqnarray*}
Clearly, $0\leq t_1\leq t_2\leq T$.
We claim that
%
\begin{equation}
\label{edensityvanish} \mbox{if $t_1<t_2$, then }
f_{\BX,T}(t)=0 \mbox{ for all $t_1<t<t_2$.}
\end{equation}
We start with the case $0<t_1<t_2<T$. Notice that, in this case,
\[
\min \bigl( f_{\BX,T}(t_1), f_{\BX,T}(t_1-)
\bigr) = \min \bigl( f_{\BX,T}(t_2), f_{\BX,T}(t_2-)
\bigr) = 0 .
\]
By \eqref{eTVaway} the density is constant on the interval
$(t_1,t_2)$. If $f_{\BX,T}(t_1)=0$, then by the right continuity of
the density, the constant must be equal to zero, so~\eqref{edensityvanish} is immediate. If $f_{\BX,T}(t_1-)=0$, then
given $\vep>0$, choose $0<s<t_1$ such that $f_{\BX,T}(s)\leq\vep$. By
\eqref{eTVaway} we know that $\mathrm{TV}_{(s, t_2)}(f_{\BX,T}) \leq
\vep$, which implies that $f(t)\leq2\vep$ on $(s, t_2)$, hence
also on $(t_1,t_2)$. Letting $\vep\to0$ proves
\eqref{edensityvanish}. If either $t_1=0$ and/or $t_2=T$, then
\eqref{edensityvanish} can be
proved using a similar argument, and the continuity of the density at 0
and at $T$ shown in part (a) of the theorem. Furthermore, we also have
%
\begin{equation}
\label{edensityvanish1} \mbox{if $t_1=t_2$, then } \min
\bigl( f_{\BX,T}(t_1), f_{\BX
,T}(t_1-)
\bigr)=0 ,
\end{equation}
with the obvious conventions in the case $t_1=t_2$ coincide with one
of the endpoints of the interval.

It follows from \eqref{edensityvanish}, \eqref{edensityvanish1}
and Lemma~\ref{lkey} that for any $\Delta>0$,
%
\begin{equation}
\label{edensvanishbig} f_{\BX,T+\Delta}(t)=0\qquad \mbox{for $t_1<t<t_2+
\Delta$.}
\end{equation}
Furthermore, we know by Lemma~\ref{llocationvar} that
%
\begin{equation}
\label{eleftcomp} F_{\BX,T+\Delta}\bigl([0,t_1]\bigr)\leq
F_{\BX,T}\bigl([0,t_1]\bigr)
\end{equation}
and
%
\begin{equation}
\label{erightcomp} F_{\BX,T+\Delta}\bigl([t_2+\Delta,T+\Delta]
\bigr)\leq F_{\BX,T}\bigl([t_2,T]\bigr) .
\end{equation}
Note that for $\Delta>0$ all the quantities in the above equations
refer to the leftmost location $\tau_{\BX,T+\Delta}$ of the supremum,
which is no longer assumed to be unique.

Since the distributions $F_{\BX,T}$ and $F_{\BX,T+\Delta}$ have equal
total masses (equal to one), it follows from
\eqref{edensvanishbig}, \eqref{eleftcomp} and
\eqref{erightcomp} that the latter two inequalities must hold as
equalities for all relevant sets. We concentrate on the resulting
equation
%
\begin{equation}
\label{eequality} F_{\BX,T+\Delta}\bigl([t_2+\Delta,T+\Delta]
\bigr)= F_{\BX,T}\bigl([t_2,T]\bigr) .
\end{equation}
Since we are working with the leftmost supremum location on a larger
interval, we can write for $\Delta>0$
\begin{eqnarray*}
P \bigl(\tau_{\BX,T}\in[t_2,T] \bigr)&=& P \bigl(
\tau_{\BX,[-\Delta,T]} \in[t_2,T] \bigr)
\\
&&{}+ P \bigl(\tau_{\BX,T}\in[t_2,T], \tau_{\BX,[-\Delta,T]} \in [-
\Delta,0) \bigr).
\end{eqnarray*}
Using Lemma~\ref{llocationvar} and \eqref{eequality}, we see that
\[
P \bigl(\tau_{\BX,T}\in[t_2,T], \tau_{\BX,[-\Delta,T]} \in [-
\Delta,0) \bigr) = 0 ,
\]
which implies that if $\Delta>T-t_2$, then
%
\begin{equation}
\label{enotlargerleft} P \Bigl(\tau_{\BX,T}\in[t_2,T],
\sup_{-\Delta\leq t\leq
-\Delta+T-t_2}X(t)\geq\sup_{t_2\leq t\leq T}X(t) \Bigr)=0 .
\end{equation}
Pick $\delta>T$. Using \eqref{enotlargerleft} with
$\Delta=n\delta-t_2, n=1,2,\ldots,$ we see that
\[
Y_n<Y_0 \qquad\mbox{a.e. on $\bigl\{ \tau_{\BX,T}
\in[t_2,T]\bigr\}$ for $n=1,2,\ldots,$}
\]
where $Y_n = \sup_{t_2-n\delta\leq t\leq T-n\delta}X(t),
n=0,1,2,\ldots.$ Note, however, that the sequence $(Y_n,
n=0,1,2,\ldots)$ is stationary, and for a stationary sequence it is
impossible that, on a set of positive probability, $Y_0>Y_n$ for
$n=1,2,\ldots$ (this is clear for an ergodic sequence; in general one
can use the ergodic decomposition). We conclude that
%
\begin{equation}
\label{enoleft} P \bigl(\tau_{\BX,T}\in[t_2,T] \bigr)=0 .
\end{equation}
Reversing the direction of time (or, equivalently, switching to the
rightmost supremum location on a larger
interval) and using Assumption~\ref{assu}, we also have
%
\begin{equation}
\label{enoright} P \bigl(\tau_{\BX,T}\in[0,t_1] \bigr)=0 .
\end{equation}

However, \eqref{edensityvanish}, \eqref{enoleft} and
\eqref{enoright} rule out any possible mass of the distribution
$F_{\BX,T}$. This contradiction shows that, under Assumption~\ref{assu},
the version of the density given by the right derivative of the cdf
$F_{\BX,T}$ is bounded away from zero. This completes the proof of the
theorem.
\end{pf}

\begin{remark} \label{rkzerodensity}
The following example shows that the statement of part (c) of Theorem
\ref{tdensityrcll} may fail without Assumption~\ref{assu}.

Let $(x(t), t\in\bbr)$ be a continuous periodic function with period
$1$, for
which $t=0$ is a global maximum. Let $U$ be a standard uniform random
variable. Then $(X(t) = x(t+U),
t\in\bbr)$ is a continuous stationary process, that always attains its
global maximum in the interval $[0,1]$. Therefore, with $T>1$, we
have $f_{\BX,T}(t)=0$ for $1\leq t<T$.
\end{remark}

Next we describe what extra restrictions on the distribution of the
location of the
supremum, in addition to the statements of Theorem
\ref{tdensityrcll}, Assumption~\ref{assl} of Section~\ref{secassumptions} imposes.
Again, one of the statements of
the theorem requires Assumption~\ref{assu}. See Remark
\ref{rknomass} for a discussion.

\begin{theorem} \label{tdensityL}
Let $\BX= (X(t), t\in\bbr)$ be a stationary sample upper
semi-continuous process, satisfying Assumption~\ref{assl}. Then the version of
the density
$f_{\BX,T}$ of the leftmost location of the supremum in the interval $[0,T]$
described in Theorem~\ref{tdensityrcll} has the following
additional properties:
\begin{longlist}[(a)]
\item[(a)] $f_{\BX,T}(0+)<\infty$, $f_{\BX,T}(T-)<\infty$ and $\mathrm{TV}_{(0,
T)}(f_{\BX,T}) \leq f_{\BX,T}(0+) +\break  f_{\BX,T}(T-)$. In particular,
the density has a bounded variation on the entire interval $(0,T)$.

\item[(b)] Assume additionally that the process is sample continuous and
satisfies Assumption~\ref{assu}.
Then either $f_{\BX,T}(t)=1/T$ for all $0<t<T$, or $\int_0^T
f_{\BX,T}(t) \,dt<1$.
\end{longlist}
\end{theorem}

\begin{remark}
Theorem~\ref{tdensityL} provides a list of specific
conditions that the
distribution of the supremum location has to satisfy (under
Assumptions~\ref{assu} and~\ref{assl}). The list turns out to be complete. That is,
for any function $f$
satisfying the conditions described in the theorem, there is a
sample continuous stationary process satisfying Assumptions~\ref{assu} and~\ref{assl},
for which $f$ is the density of the supremum
location. Thus we have obtained a full characterization of the set
of all possible densities. In order to decide whether a candidate
function can be the density of the supremum location for some
stationary process, we only need to check the list of
conditions given in the theorem. This is, of course a much
easier task than trying to construct an appropriate process. We refer
the reader to
\citet{samorodnitskyshen2012} for details and proofs.
\end{remark}

\begin{remark}
Note that part (b) of Theorem~\ref{tdensityL} says that, unless
the location of the supremum is uniformly distributed in the
interval $(0,T)$, the supremum is achieved, with a positive
probability, at an endpoint of the interval. The proof of this
part, exhibited in the following pages, actually implies more. It
shows that the uniform distribution occurs only
when the suprema of the process appear periodically with period
equal to $T$:
\[
P \bigl( X ( \tau_{\BX,[T,2T]} )= X ( \tau_{\BX,T} ),
\tau_{\BX,[T,2T]}-\tau_{\BX,T}=T \bigr)=1.
\]
\end{remark}

\begin{pf*}{Proof of Theorem~\ref{tdensityL}}
Assumption~\ref{assl} and stationarity imply that for any $0<t<T$,
\begin{eqnarray*}
f_{\BX,T}(t) &=& \lim_{\vep\downarrow0}\frac{P(\tau_{\BX,T}\in
(t,t+\vep))}{\vep}
\\
&\leq&\limsup_{\vep\downarrow0}\frac{P( \mbox{$\BX$ has a local
maximum in $(t,t+\vep)$)}}{\vep}
\\
&=& \limsup_{\vep\downarrow0}\frac{P( \mbox{$\BX$ has a local
maximum in $(0,\vep)$)}}{\vep}\leq K .
\end{eqnarray*}
This proves finiteness of $f_{\BX,T}(0+)<\infty$ and
$f_{\BX,T}(T-)$. The rest of the statement in part (a) follows from
\eqref{eTVleft1} by letting $\vep\uparrow T$.

We now prove part (b). Assume that $P(\tau_{\BX,T}=0 \mbox{ or } T)=0$.
By stationarity this implies that $\tau_{\BX,[T,2T]}\in(T,2T)$
with probability 1. We first prove that
%
\begin{equation}
\label{esmallnoteq} P \bigl( X ( \tau_{\BX,[T,2T]} )\not= X (
\tau_{\BX,T} ) \bigr)=0 .
\end{equation}
By symmetry, it is enough to prove the one-sided claim
%
\begin{equation}
\label{esmallnoteq1} P \bigl( X ( \tau_{\BX,[T,2T]} )<X ( \tau_{\BX,T}
) \bigr)=0 .
\end{equation}
Indeed, suppose, to the contrary, that the probability in
\eqref{esmallnoteq1} is positive. Under Assumption~\ref{assu} we can use
the continuity from below of measures to see that there is $\vep>0$
such that
\[
p:=P \Bigl( X ( \tau_{\BX,T} ) >X ( \tau_{\BX,[T,2T]} )+\vep, X (
\tau_{\BX,T} ) > \max_{t\in L_T, t\not= \tau_{\BX,T}}X(t)+\vep \Bigr)>0 .
\]
Here $L_T$ is the (a.s. finite) set of the local maxima of $\BX$ in
the interval $(0,T)$.

Next, by the uniform continuity of the process $\BX$ on $[0,T]$, there
is $n\geq1$ such that
\[
P \Bigl( \sup_{0\leq s<t\leq T, t-s\leq T/n} \bigl| X(t)-X(s) \bigr|>\vep/2 \Bigr)\leq p/2 .\vadjust{\goodbreak}
\]
We immediately conclude by the law of total probability that there is
$i=1,\ldots,n$ such that $P(A_i)>0$, where
\begin{eqnarray*}
A_i &=& \bigl\{ X ( \tau_{\BX,T} ) >X ( \tau_{\BX,[T,2T]} )+
\vep, X ( \tau_{\BX,T} ) > \max_{t\in L_T, t\not= \tau_{\BX,T}}X(t)+\vep,
\\
&&\hspace*{5pt}{}(i-1)T/n<\tau_{\BX,T}<iT/n, \sup_{(i-1)T/n\leq s,t\leq iT/n} \bigl| X(t)-X(s) \bigr|\leq\vep/2
\bigr\}.
\end{eqnarray*}
However, on the event $A_i$, $X(iT/n)= \sup_{iT/n\leq t\leq2T}X(t)$,
implying that $\tau_{\BX, [iT/n,iT/n+T]}= iT/n$. By stationarity, this
contradicts the assumption\break $P(\tau_{\BX,T}=0)=0$. This contradiction
proves \eqref{esmallnoteq1} and, hence, also
\eqref{esmallnoteq}.

Next, we check that
%
\begin{equation}
\label{esmalleq} P \bigl( X ( \tau_{\BX,[T,2T]} )= X ( \tau_{\BX,T} ),
\tau_{\BX,[T,2T]}-\tau_{\BX,T}<T \bigr)=0 .
\end{equation}
Indeed, suppose that, to the contrary, the probability above is
positive. By the continuity from below of measures, there is
$\vep>0$ such that
\[
P \bigl( X ( \tau_{\BX,[T,2T]} )= X ( \tau_{\BX,T} ),
\tau_{\BX,[T,2T]}-\tau_{\BX,T}<T-\vep \bigr)>0 .
\]
Take $n>2T/\vep$. By the law of total probability there are
$i_1,i_2=1,\ldots,n$ such that $P(A_{i_1,i_2})>0$, where
\begin{eqnarray*}
A_{i_1,i_2} &=& \bigl\{ X ( \tau_{\BX,[T,2T]} )= X ( \tau_{\BX,T} ),
\tau_{\BX,[T,2T]}-\tau_{\BX,T}<T-\vep,
\\
&&\hspace*{5pt}{}(i_1-1)T/n<\tau_{\BX,T}<i_1T/n,\\
&&\hspace*{5pt}{}T+(i_2-1)T/n<\tau_{\BX,[T,2T]}<T+i_2T/n \bigr\}.
\end{eqnarray*}
By the choice of $n$, $T+i_2T/n -(i_1-1)T/n<T$, so that, on the event
$A_{i_1,i_2}$, the process $\BX$ has at least two points,
$\tau_{\BX,T}$ and $\tau_{\BX,[T,2T]}$, at which the supremum over the
interval $[(i_1-1)T/n, (i_1-1)T/n+T]$ is achieved. By stationarity,
this contradicts Assumption~\ref{assu}. This contradiction proves~\eqref{esmalleq}.

Finally, we check that
%
\begin{equation}
\label{elargeeq} P \bigl( X ( \tau_{\BX,[T,2T]} )= X ( \tau_{\BX,T} ),
\tau_{\BX,[T,2T]}-\tau_{\BX,T}>T \bigr)=0 .
\end{equation}
The proof is similar to the proof of \eqref{esmallnoteq1}, so we
only sketch the argument. Suppose that, to the contrary, the
probability in
\eqref{elargeeq} is positive. Use the continuity of measures to see
that the probability remains positive if we require that
$\tau_{\BX,[T,2T]}-\tau_{\BX,T}>T+\vep$ for some $\vep>0$. Next, use
Assumption~\ref{assu} to separate the value of $X ( \tau_{\BX,T} )$
from the values of $\BX$ at other local maxima in $(0,T)$ and,
finally, use the uniform continuity of the process $\BX$ to show that
there is a point $T<b<2T$ and an event of positive probability on
which $\tau_{\BX,[b-T,b]}=b$. By stationarity, this contradicts the
assumption $P(\tau_{\BX,T}=T)=0$.

Combining \eqref{esmallnoteq}, \eqref{esmalleq} and
\eqref{elargeeq}, we see that the assumption $P(\tau_{\BX,T}=0 \mbox{ or } T)=0$ implies that
%
\begin{equation}
\label{eallequal} P \bigl( X ( \tau_{\BX,[T,2T]} )= X ( \tau_{\BX,T}
), \tau_{\BX,[T,2T]}-\tau_{\BX,T}=T \bigr)=1.
\end{equation}
Let $0<a<b<T$. We have by stationarity,
\begin{eqnarray*}
P \bigl( \tau_{\BX,T}\in(0,b-a) \bigr) &=& P \bigl( \tau_{\BX,[a,a+T]}
\in(a,b) \bigr)
\\
&=& P \bigl( \tau_{\BX,[a,a+T]}\in(a,b), \tau_{\BX,T}\in(0,a) \bigr) \\
&&{}+ P
\bigl( \tau_{\BX,[a,a+T]}\in(a,b), \tau_{\BX,T}\in(a,T) \bigr) .
\end{eqnarray*}
By \eqref{eallequal}, if $\tau_{\BX,T}\in(0,a)$, then
$\tau_{\BX,[T,2T]}\in(T,T+a)$ and $X (
\tau_{\BX,[T,2T]} )>\sup_{t\in[a,b]} X(t)$. Therefore, the first
term in the right-hand side above vanishes. Similarly, by
\eqref{eallequal}, if $\tau_{\BX,T}\in(a,T)$, then
$\tau_{\BX,[T,2T]}\in(T+a,2T)$, and $X ( \tau_{\BX,T}
)>\sup_{t\in
[T,T+a]}X(t)$. Therefore,
\[
P \bigl( \tau_{\BX,T}\in(0,b-a) \bigr) = P \bigl( \tau_{\BX,T}
\in(a,b) \bigr)
\]
for any $0<a<b<T$, which proves the uniformity of the distribution of~$\tau_{\BX,T}$.
\end{pf*}

\begin{remark} \label{rknomass}
A simple special case of the process in Remark~\ref{rkzerodensity}
shows that the statement of part (b) of Theorem~\ref{tdensityL} may
fail without Assumption~\ref{assu}.

We take, for clarity, a specific function $x$. Let $x(t) = 1-2|t|$ for
$|t|\leq1/2$ and extend $x$ to a periodic function with period
$1$. Then for any $T>1$, the leftmost location of the supremum in the
interval $[0,T]$ of the
process $(X(t) = x(t+U), t\in\bbr)$ is in the interval $(0,1)$ with
probability 1, and (as we already know) this location is not uniformly
distributed between 0 and~$T$.
\end{remark}

None of the statement of Theorem~\ref{tdensityL} holds, in general,
without Assumption~\ref{assl}, as the following example shows.

\begin{example} \label{exO-U}
Let $X(t) = e^{-t/2}B(e^t), t\geq0$, where $(B(t))$ is the standard
Brownian motion. Then $\BX$ is a stationary Gaussian process, the
Ornstein--Uhlenbeck process. It is, clearly, sample continuous, and the
strong Markov property of the Brownian motion shows that, for any
$T>0$, it satisfies Assumption~\ref{assu}. It is clear that Assumption~\ref{assl}
fails for the Ornstein--Uhlenbeck process.

By the law of iterated logarithm for the Brownian motion we see that,
on a set of probability 1, in any interval $(0,\vep)$ with $\vep>0$
there is a point $t$ such that $X(t)>X(0)$. Therefore,
$P(\tau_{\BX,T}=0)=0$ and, similarly, $P(\tau_{\BX,T}=T)=0$ for any
$T>0$.

It is also easy to show, using the basic properties of the Brownian
motion, that the density $f_{\BX,T}$ is not bounded near each of the
two endpoints of the interval $[0,T]$, so that both statements of
Theorem~\ref{tdensityL} fail for this process.
\end{example}

\section{Universal upper bounds on the density}
\label{secupperbounds}

The upper bounds in part (b) of Theorem~\ref{tdensityrcll} turn out to
be the best
possible pointwise, as is shown in the following result.\vadjust{\goodbreak}
%
\begin{proposition} \label{pruppergen}
For each $0<t<T$ and any number
smaller than the upper bound given in \eqref{edensitybound},
there is a sample continuous stationary process satisfying Assumptions~\ref{assu} and~\ref{assl} for which the right continuous version of the
density $f_{\BX,T}(t)$ of the supremum location at time $t$ exceeds
that number.
\end{proposition}
\begin{pf}
By symmetry, it is enough to show that for
any $0<t<T$ and any number smaller than $1/t$ there is a stationary
process of the required type for which $f_{\BX,T}(t)$ exceeds that
number.

To this end, let $\tau>t$ and let $k\geq1$ be an integer.
We define a periodic function $(x(s), s\in\bbr)$
with period $k\tau+2T$ by defining its values on the interval
$[0,k\tau+2T]$. We set $x(i\tau)= k-i$ for $i=0,1,\ldots, k$ and
$x(k\tau+2T)= k$. We set, further, for $i=0,1,\ldots, k-1$,
$x((i+1/2)\tau)=-R$ and also $x ( k\tau+T )=-R$
for a large positive $R$ we describe in a moment. We complete
the definition of the function by connecting linearly the values in
neighboring points where the function has already been
defined. Fix $t<r<\tau$, and choose now $R$ so large that the condition
%
\begin{equation}
\label{elocalmax} x ( i\tau )>x ( i\tau-r )
\end{equation}
holds for all $i=1,\ldots, k$. Now define a stationary process by
$X(s) = x(s-U), s\in\bbr$, where $U$ is uniformly distributed
between 0 and $k\tau+ 2T$. By construction, the process is sample
continuous and satisfies Assumptions~\ref{assu} and~\ref{assl}.

If, for $i=1,\ldots, k$, we have $i\tau-r<U<i\tau$, then the local
maximum at $s=i\tau$ of the function $\mathbf{x}$ becomes the global maximum
of the process $\BX$ over the interval $[0,T]$, and is located in
the interval $(0,r)$. This contributes $1/(k\tau+2T)$ to the value of
the density $f_{\BX,T} $ at each point of the interval $(0,r)$. In
particular, since $t\in(0,r)$,
\[
f_{\BX,T}(t)\geq\frac{k}{k\tau+2T} .
\]
Since we can take $k$ arbitrarily large, the value of the density can
be arbitrarily close to $1/\tau$, and since $\tau$ can be taken
arbitrarily close to $t$, the value of the density can be arbitrarily
close to $1/t$.
\end{pf}

Suppose now that the stationary process $\BX$ is time
reversible, that is, if $(X(-t), t\in\bbr)\eid(X(t),
t\in\bbr)$. That would, obviously, be the case for stationary Gaussian
processes. If the process satisfies also Assumption~\ref{assu}, then the
distribution of the unique supremum location $\tau_{\BX, T}$ is
symmetric in the interval $[0,T]$, meaning that $\tau_{\BX, T}\eid
T-\tau_{\BX, T}$. Therefore, the density $f_{\BX, T}$ satisfies
%
\begin{equation}
\label{esymmdensity} f_{\BX, T}(t) = f_{\BX, T}(T-t)
\end{equation}
for all $0<t<T/2$ that are continuity points of $f_{\BX, T}$.
Even though the upper bound given in part (b) of Theorem
\ref{tdensityrcll} is symmetric around the middle of the interval
$[0,T]$, it\vadjust{\goodbreak} turns out that the bounded variation property in part (d)
of Theorem
\ref{tdensityrcll} provides a better bound in this symmetric
case. This bound and its optimality, even within the class of
stationary Gaussian processes, is presented in the following result.
%
\begin{proposition}\label{csymmdensitybound}
Let $\BX= (X(t), t\in\bbr)$ be a time reversible stationary sample upper
semi-continuous process satisfying Assumption~\ref{assu}. Then the density
$f_{\BX,T}$ of the unique location of the supremum in the interval
$[0,T]$ satisfies
%
\begin{equation}
\label{esymmdensitybound} f_{\BX,T}(t) \leq\cases{ %
\displaystyle\frac{1}{2t}, & \quad $\mbox{if $\displaystyle 0<t\leq\frac{T}{3}$},$
\vspace*{2pt}\cr
\displaystyle\frac{1}{T-t}, & \quad $\mbox {if $\displaystyle \frac{T}{3}<t\leq\frac{T}{2}$},$
\vspace*{2pt}\cr
\displaystyle\frac1t, & \quad $\mbox{if $\displaystyle \frac{T}{2}<t\leq\frac{2T}{3}$},$
\vspace*{2pt}\cr
\displaystyle\frac{1}{2(T-t)}, &\quad $\mbox{if $\displaystyle \frac{2T}{3}<t<T$.}$}
\end{equation}

Furthermore, for each $0<t<T$ and any number
smaller than the upper bound given in \eqref{esymmdensitybound},
there is a sample continuous Gaussian
process satisfying Assumptions~\ref{assu} and~\ref{assl} for
which the density $f_{\BX,T}(t)$ exceeds that number.
\end{proposition}
\begin{pf}
Since the density $f_{\BX,T}$ is right continuous, it is enough to
consider only continuity points of the density, and by
\eqref{esymmdensity}, it is enough to consider $0<t<T/2$. Then $T-t$
is also a continuity point of the density. Denote
$a=\inf_{0<s\leq t}f_{\BX,T}(s)$, $b=\inf_{t<s<T/2}f_{\BX,T}(s)$. Note
that, given $\vep>0$, there is a continuity point of the density $u\in
(0,t]$ such that
$f_{\BX,T}(u)\leq a+\vep$, and there is a continuity point of the
density $v\in[t,T/2]$ such
that $f_{\BX,T}(v)\leq b+\vep$. Observe also that
%
\begin{equation}
\label{eboundary} at + b(T/2-t) \leq\int_0^{T/2}
f_{\BX,T}(s) \,ds \leq\frac12 .
\end{equation}
Furthermore, applying the total variation bound \eqref{eTVaway} to
the interval $[u,\break T-u]$ gives us
\begin{eqnarray*}
2(a+\vep) &\geq& f_{\BX,T}(u) + f_{\BX,T}(T-u)
\\
&\geq& \bigl| f_{\BX,T}(t)- f_{\BX,T}(u) \bigr| + \bigl| f_{\BX,T}(v)-
f_{\BX,T}(t) \bigr|
\\
&&{}+ \bigl| f_{\BX,T}(T-v)- f_{\BX,T}(v)\bigr | +\bigl | f_{\BX,T}(T-t)-
f_{\BX,T}(T-v) \bigr|
\\
&&{}+ \bigl| f_{\BX,T}(T-u)- f_{\BX,T}(T-t) \bigr|
\\
&\geq& 2 \bigl( f_{\BX,T}(t)-a-\vep \bigr)_+ + 2 \bigl( f_{\BX,T}(t)-b-
\vep \bigr)_+ .
\end{eqnarray*}
Letting $\vep\to0$ and recalling that $a\leq f_{\BX,T}(t)$ and
$b\leq
f_{\BX,T}(t)$, we obtain
%
\begin{equation}
\label{etargetfunction} f_{\BX,T}(t) \leq a+b/2 .\vadjust{\goodbreak}
\end{equation}
Since $b\leq f_{\BX,T}(t)$, this implies that
%
\begin{equation}
\label{eboundary1} b\leq2a .
\end{equation}
If $0<t\leq T/3$, then the largest value of the right-hand side of
\eqref{etargetfunction} under the constraint \eqref{eboundary}
requires taking $a$ as large as possible and $b$ as small as
possible. Taking $a=1/2t$ and $b=0$ in \eqref{etargetfunction}
results in the upper bound given in \eqref{esymmdensitybound} in
this range. If $T/3<t\leq T/2$, then the largest value of the
right-hand side of
\eqref{etargetfunction} under the constraint \eqref{eboundary}
requires taking $a$ as small as possible and $b$ as large as
possible. By \eqref{eboundary1}, we have to take $a=1/2(T-t)$,
$b=1/(T-t)$ in \eqref{etargetfunction}, which results in the upper
bound given in \eqref{esymmdensitybound} in this case.

It remains to prove the optimality part of the statement of the
corollary. By symmetry it is enough to consider $0<t\leq T/2$. Fix
such $t$. Let $\vep>0$ be a small number and $h>0$ be a large number,
rationally independent of $t+\vep$. Consider a stationary Gaussian process
given by
\begin{eqnarray*}
X(s) &=& G_1\cos \biggl( \frac{2\pi}{t+\vep}s \biggr) +
G_2\sin \biggl( \frac{2\pi}{t+\vep}s \biggr)
\\
&&{}+ G_3\cos \biggl( \frac{2\pi}{h}s \biggr) + G_4
\sin \biggl( \frac{2\pi}{h}s \biggr), \qquad s\in\bbr,
\end{eqnarray*}
where $G_1,\ldots, G_4$ are i.i.d. standard normal random
variables. The process is, clearly, sample continuous, and it
satisfies Assumption~\ref{assl}. Furthermore, rational independence of $t+\vep$
and $h$ implies that, on a set of probability 1, the process $\BX$ has
different values at all of its local maxima, hence Assumption~\ref{assu} is
satisfied for any $T>0$. Note that we can write
\[
X(s) = A_1\cos \biggl( \frac{2\pi}{t+\vep}s +U_1
\biggr) + A_2\cos \biggl( \frac{2\pi}{h}s+U_2 \biggr)
:= X_1(s) + X_2(s),\qquad s\in\bbr,
\]
where $A_1$ and $A_2$ have the density $xe^{-x^2/2}$ on $(0,\infty)$,
and $U_1$ and $U_2$ are uniformly distributed between 0 and $2\pi$,
with all 4 random variables being independent. Clearly, the leftmost
location of the supremum of the process $\BX_1$ is at
\[
\tau_1= (t+\vep)\frac{2\pi-U_1}{2\pi},
\]
which is uniformly distributed between 0 and $t+\vep$. On the event
$E=\{ 0<U_2<\pi-2\pi T/h\}$ the process $\BX_2$ is decreasing on
$[0,T]$, so the value of the sum $\BX$ at the leftmost supremum of
$\BX_1$ exceeds the value of the sum at all the other locations of the
supremum of $\BX_1$ in the interval $[0,T]$. If the supremum of the
sum remained at $\tau_1$, the density of that unique supremum would be
at least $P(E)/(t+\vep)$ at each point of the interval
$(0,t+\vep)$. Since $P(E)\to1/2$ as $h\to\infty$, the value of the
density at $t$ would exceed any value smaller than $1/2t$ after taking
$h$ large and $\vep$ small. The location of the supremum of the sum
does not remain at $\tau_1$ but, instead, moves to
$\tau_2=\tau_2(A_1,A_2,U_1,U_2)$ defined by
\[
\tau_2 = \sup \biggl\{s\leq\tau_1\dvtx
\frac{A_1}{t+\vep}\sin \biggl( \frac
{2\pi}{t+\vep}s +U_1 \biggr) +
\frac{A_2}{h}\sin \biggl( \frac{2\pi}{h}s+U_2 \biggr)=0
\biggr\}.
\]
For large $h$, $\tau_2$ is nearly identical to $\tau_1$, and
straightforward but somewhat tedious calculus based on the implicit
function theorem shows that the above statement remains true for
$\tau_2$: the contribution of the event $E$ to the density of the
unique supremum of the process $\BX$ would exceed any value smaller
than $1/2t$ at any point of the interval $(0,t+\vep)$ after taking
$h$ large and $\vep$ small. We omit the details.

We have shown the optimality of the upper bound given in
\eqref{esymmdensitybound} in the case $0<t\leq T/3$.
It remains to consider the case $T/3<t\leq T/2$. We will use again a
two-wave stationary Gaussian process, but with a slightly different
twist. Let $\vep>0$ be a small number, $h>0$ a large number and $r>0$
a fixed number that is rationally independent of $T-t+\vep$. Consider a
stationary Gaussian process given by
\begin{eqnarray*}
X(s) &=& A_1\cos \biggl( \frac{2\pi}{T-t+\vep}s+U_1 \biggr)
+\frac1h A_2\cos \biggl( \frac{2\pi}{r}s+U_2 \biggr)
\\
&:=& X_1(s)+X_2(s), s\in\bbr,
\end{eqnarray*}
where $A_1,A_2,U_1$ and $U_2$ are as above. As above, $\BX$ is a
sample continuous Gaussian process satisfying Assumptions~\ref{assl} and~\ref{assu}. Now the leftmost location
of the supremum of the process $\BX_1$ is at
\[
\tau_1 = (T-t+\vep)\frac{2\pi-U_1}{2\pi} ,
\]
which is uniformly distributed between 0 and $T-t+\vep$. Further, if
$\tau_1>t-\vep/2$, then $\tau_1$ is the unique supremum of $\BX_1$ in
the interval $[0,T]$. If the supremum of the sum $\BX$ remained at
$\tau_1$, then the density of the supremum location at the point $t$
would be at least
$1/(T-t+\vep)$, which would then exceed any value smaller than
$1/(T-t)$ after taking $\vep$ small. The location of the supremum of
$\BX$ does not remain at $\tau_1$, but instead moves to the unique for
large $h$ point $\tau_2=\tau_2(A_1,A_2,U_1,U_2)$ in $[0,T]$ satisfying
\[
\frac{A_1}{T-t+\vep}\sin \biggl( \frac{2\pi}{T-t+\vep}\tau_2+U_1
\biggr) + \frac{A_2}{hr}\sin \biggl( \frac{2\pi}{r}\tau_2+U_2
\biggr)=0 .
\]
For large $h$, $\tau_2$ is nearly identical to $\tau_1$, and as above,
using the implicit value theorem allows us to conclude that, for any
value smaller than $1/(T-t)$, the value of the density of $\tau_2$ in
the interval $(t-\vep/2,T-t+\vep)$ exceeds that value after taking
$\vep$ small and $h$ large. This proves the optimality of the upper
bound given in~\eqref{esymmdensitybound} in all cases.
\end{pf}

%


\printaddresses


\begin{thebibliography}{7}

\bibitem[\protect\citeauthoryear{Adler}{1990}]{adler1990}
\begin{bbook}[mr]
\bauthor{\bsnm{Adler},~\bfnm{Robert~J.}\binits{R.~J.}}
(\byear{1990}).
\btitle{An Introduction to Continuity, Extrema, and Related Topics for General
  {G}aussian Processes}.
\bseries{Institute of Mathematical Statistics Lecture Notes---Monograph Series}
\bvolume{12}.
\bpublisher{IMS}, \blocation{Hayward, CA}.
\bid{mr={1088478}}
\bptok{imsref}%
\end{bbook}
\endbibitem

\bibitem[\protect\citeauthoryear{Adler and Taylor}{2007}]{adlertaylor2007}
\begin{bbook}[mr]
\bauthor{\bsnm{Adler},~\bfnm{Robert~J.}\binits{R.~J.}} \AND
  \bauthor{\bsnm{Taylor},~\bfnm{Jonathan~E.}\binits{J.~E.}}
(\byear{2007}).
\btitle{Random Fields and Geometry}.
\bpublisher{Springer}, \blocation{New York}.
\bid{mr={2319516}}
\bptok{imsref}%
\end{bbook}
\endbibitem

\bibitem[\protect\citeauthoryear{Aza{\"{\i}}s and
  Wschebor}{2002}]{azaiswschebor2002}
\begin{bincollection}[mr]
\bauthor{\bsnm{Aza{\"{\i}}s},~\bfnm{Jean-Marc}\binits{J.-M.}} \AND
  \bauthor{\bsnm{Wschebor},~\bfnm{Mario}\binits{M.}}
(\byear{2002}).
\btitle{The distribution of the maximum of a {G}aussian process: {R}ice method
  revisited}.
In \bbooktitle{In and Out of Equilibrium ({M}ambucaba, 2000)}.
\bseries{Progress in Probability}
\bvolume{51}
\bpages{321--348}.
\bpublisher{Birkh\"auser}, \blocation{Boston, MA}.
\bid{mr={1901961}}
\bptok{imsref}%
\end{bincollection}
\endbibitem

\bibitem[\protect\citeauthoryear{Aza{\"{\i}}s and
  Wschebor}{2009}]{azaiswschebor2009}
\begin{bbook}[mr]
\bauthor{\bsnm{Aza{\"{\i}}s},~\bfnm{Jean-Marc}\binits{J.-M.}} \AND
  \bauthor{\bsnm{Wschebor},~\bfnm{Mario}\binits{M.}}
(\byear{2009}).
\btitle{Level Sets and Extrema of Random Processes and Fields}.
\bpublisher{Wiley}, \blocation{Hoboken, NJ}.
\bid{doi={10.1002/9780470434642}, mr={2478201}}
\bptok{imsref}%
\end{bbook}
\endbibitem

\bibitem[\protect\citeauthoryear{Leadbetter, Lindgren and
  Rootz{\'e}n}{1983}]{leadbetterlindgrenrootzen1983}
\begin{bbook}[mr]
\bauthor{\bsnm{Leadbetter},~\bfnm{M.~R.}\binits{M.~R.}},
  \bauthor{\bsnm{Lindgren},~\bfnm{Georg}\binits{G.}} \AND
  \bauthor{\bsnm{Rootz{\'e}n},~\bfnm{Holger}\binits{H.}}
(\byear{1983}).
\btitle{Extremes and Related Properties of Random Sequences and Processes}.
\bpublisher{Springer}, \blocation{New York}.
\bid{mr={0691492}}
\bptok{imsref}%
\end{bbook}
\endbibitem

\bibitem[\protect\citeauthoryear{Samorodnitsky and Shen}{2012}]{samorodnitskyshen2012}
\begin{barticle}[mr]
\bauthor{\bsnm{Samorodnitsky},~\bfnm{Gennady}\binits{G.}}
\AND
\bauthor{\bsnm{Shen},~\bfnm{Yi}\binits{Y.}}
(\byear{2012}).
\btitle{Distribution of the supremum location of stationary processes}.
\bjournal{Electron. J. Probab.}
\bvolume{17}
\bpages{1--17}.
\bptok{imsref}%
\end{barticle}
\endbibitem

\end{thebibliography}
\end{document}